\documentclass[12pt]{article} 
\usepackage{amssymb,latexsym,color,psfig}
\newcounter{environment}[section]
\renewcommand{\theenvironment}{%
\arabic{section}.\arabic{environment}}

{\begin{rm}\refstepcounter{environment}{\textbf\theenvironment\
\bf Definition.~~}}%
{\end{rm}}

{\begin{rm}\refstepcounter{environment}{\textbf\theenvironment\
\bf Example.~~}}%
{\end{rm}}

\newenvironment{conjecture}%
{\begin{rm}\refstepcounter{environment}{\textbf\theenvironment\
\bf Conjecture.~~}}%
{\end{rm}}

\newenvironment{proposition}%
{\begin{rm}\refstepcounter{environment}{\textbf\theenvironment\
\bf Proposition.~~}}%
{\end{rm}}

{\begin{rm}\refstepcounter{environment}{\textbf\theenvironment\
\bf Proposition and Definition.~~}}%
{\end{rm}}

\newenvironment{theorem}%
{\begin{rm}\refstepcounter{environment}{\textbf\theenvironment\
\bf Theorem.~~}}%
{\end{rm}}

\newenvironment{corollary}%
{\begin{rm}\refstepcounter{environment}{\textbf\theenvironment\
\bf Corollary.~~}}%
{\end{rm}}

\newenvironment{lemma}%
{\begin{rm}\refstepcounter{environment}{\textbf\theenvironment\
\bf Lemma.~~}}%
{\end{rm}}

\savebox0{\rule[-2\unitlength]{0pt}{10\unitlength}%
\begin{picture}(10,10)(0,2)
\put(0,0){\line(0,1){10}}
\put(0,10){\line(1,0){10}}
\put(10,0){\line(0,1){10}}
\put(0,0){\line(1,0){10}}
\end{picture}}
\savebox1{\rule[-2\unitlength]{0pt}{10\unitlength}%
\begin{picture}(10,10)(0,2)
\put(0,10){\line(1,0){10}}
\put(10,0){\line(0,1){10}}
\put(0,10){\line(1,-1){10}}
\end{picture}}

\renewcommand\emptyset{\mbox{\sc \O}}
\newlength\cellsize 
\savebox2{%
\begin{picture}(18,18)
\put(0,0){\line(1,0){18}}
\put(0,0){\line(0,1){18}}
\put(18,0){\line(0,1){18}}
\put(0,18){\line(1,0){18}}
\end{picture}}
\newcommand\cellify[1]{\def\thearg{#1}\def\nothing{}%
\ifx\thearg\nothing
\vrule width0pt height\cellsize depth0pt\else
\hbox to 0pt{\usebox2\hss}\fi%
\vbox to 18\unitlength{
\vss
\hbox to 18\unitlength{\hss$#1$\hss}
\vss}}
\newcommand\tableau[1]{\vtop{\let\\=\cr
\setlength\baselineskip{-16000pt}
\setlength\lineskiplimit{16000pt}
\setlength\lineskip{0pt}
\halign{&\cellify{##}\cr#1\crcr}}}
\savebox3{%
\begin{picture}(15,15)
\put(0,0){\line(1,0){15}}
\put(0,0){\line(0,1){15}}
\put(15,0){\line(0,1){15}}
\put(0,15){\line(1,0){15}}
\end{picture}}
\newcommand\expath[1]{%
\hbox to 0pt{\usebox3\hss}%
\vbox to 15\unitlength{
\vss
\hbox to 15\unitlength{\hss$#1$\hss}
\vss}}
\begin{document}
\newcommand{\qbc}[2]{ {\left [{#1 \atop #2}\right ]}}
\newcommand{\anbc}[2]{{\left\langle {#1 \atop #2} \right\rangle}}
\newcommand{\be}{\begin{enumerate}} \newcommand{\ee}{\end{enumerate}}
\newcommand{\beq}{\begin{equation}} \newcommand{\eeq}{\end{equation}}
\newcommand{\bea}{\begin{eqnarray}} \newcommand{\eea}{\end{eqnarray}}
\newcommand{\beas}{\begin{eqnarray*}}
  \newcommand{\eeas}{\end{eqnarray*}} \newcommand{\zz}{\mathbb{Z}}
\newcommand{\pp}{\mathbb{P}} \newcommand{\nn}{\mathbb{N}}
\newcommand{\qq}{\mathbb{Q}} \newcommand{\rr}{\mathbb{R}}
\newcommand{\bm}[1]{{\mbox{\boldmath $#1$}}}
\newcommand{\sn}{\mathfrak{S}_n} \newcommand{\cs}{{\cal S}}
\newcommand{\la}{\lambda} \newcommand{\lm}{{\lambda/\mu}}
\newcommand{\ep}{{\cal E}_P} \newcommand{\lp}{{\cal L}_P}
\newcommand{\lpw}{{\cal L}_{P,\omega}} \newcommand{\ppw}{{(P,\omega)}}
\newcommand{\dla}{\mathrm{Dom}_\lambda} \newcommand{\pw}{{P,\omega}}
\newcommand{\ds}{\displaystyle} \newcommand{\fs}{\mathfrak{S}}
\newcommand{\st}{\,:\,} \newcommand{\maj}{\mathrm{maj}}
\newcommand{\modt}{\,(\mathrm{mod}\,2)}
\newcommand{\tr}{\textcolor{red}} \newcommand{\tb}{\textcolor{blue}}
\newcommand{\tg}{\textcolor{green}}
\newcommand{\tm}{\textcolor{magenta}}
\newcommand{\tbn}{\textcolor{brown}}
\newcommand{\tp}{\textcolor{purple}}
\newcommand{\tn}{\textcolor{nice}}
\newcommand{\tor}{\textcolor{orange}}

\definecolor{brown}{cmyk}{0,0,.35,.65}
\definecolor{purple}{rgb}{.5,0,.5}
\definecolor{nice}{cmyk}{0,.5,.5,0}
\definecolor{orange}{cmyk}{0,.35,.65,0}

\begin{centering}
\textcolor{red}{\Large\bf Some Remarks
on Sign-Balanced and Maj-Balanced Posets}\\[.2in] 
\textcolor{blue}{Richard P. Stanley}\footnote{Partially supported by
  NSF grant \#DMS-9988459.}\\ 
Department of Mathematics\\
Massachusetts Institute of Technology\\
Cambridge, MA 02139\\
\emph{e-mail:} rstan@math.mit.edu\\[.2in]
\textcolor{magenta}{version of 14 January 2004}\\[.2in]

\end{centering}
\vskip 10pt
\section{Introduction.} \label{sec:intro}
Let $P$ be an $n$-element poset (partially ordered set), and let
$\omega:P\rightarrow [n]=\{1,2,\dots,n\}$ be a bijection, called a
\emph{labeling} of $P$. We call the pair $(P,\omega)$ a
\emph{labelled poset}. A \emph{linear extension} of $P$ is an order-preserving
bijection $f:P\rightarrow [n]$.  We can regard
$f$ as defining a permutation $\pi=\pi(f)$ of the set $[n]$ given by $\pi(i)=j$ if 
$f(\omega^{-1}(j))=i$. We write $\pi$ in the customary way as a word $a_1a_2\cdots
a_n$, where $\pi(i)=a_i=\omega(f^{-1}(i))$. We will say for instance
that $f$ is an \emph{even linear extension} of $(P,\omega)$ if $\pi$
is an even permutation (i.e., an 
element of the alternating group $\mathfrak{A}_n$). Let $\ep$ denote
the set of linear extensions of $P$, and set $\lpw=\{\pi(f)\st f\in\ep\}$

We say that $\ppw$ is \emph{sign-balanced} if $\lpw$ contains the same
number of even permutations as odd permutations. Note that the parity
of a linear extension $f$ depends on the labeling $\omega$.  However,
the notion of sign-balanced depends only on $P$, since changing the
labeling of $P$ simply multiplies the elements of $\lpw$ by a fixed
permutation in $\sn$, the symmetric group of all permutations of
$[n]$. Thus we can simply say that $P$ is sign-balanced without
specifying $\omega$.

We say that a function $\vartheta:\ep\rightarrow\ep$ is
\emph{parity-reversing} (respectively, \emph{parity-preserving}) if
for all $f\in\ep$, the permutations $\pi(f)$ and $\pi(\vartheta(f))$
have opposite parity (respectively, the same parity). Note that the
properties of parity-reversing and parity-preserving do not depend on
$\omega$; indeed, $\vartheta$ is parity-reversing (respectively,
parity-preserving) if and only if for all $f\in\ep$, the permutation
$\vartheta f\circ f^{-1}\in\sn$ is odd (respectively, even),

Sign-balanced posets were first considered by Ruskey \cite{ruskey}.
He established the following result, which shows that many
combinatorially occuring classes of posets, such as geometric lattices
and Eulerian posets, are sign-balanced.

\begin{theorem} \label{thm:p-r}
\emph{Suppose $\#P\geq 2$. If every nonminimal element of the poset
$P$ is greater than at least two minimal elements, then $P$ is
sign-balanced.} 
\end{theorem}

\textbf{Proof.} Let $\pi=a_1a_2a_3\cdots a_n\in\lpw$. Let $\pi'=
\pi(1,2) = a_2a_1a_3\cdots a_n\in \sn$. (We always multiply
permutations from right to left.)  By the hypothesis on $P$, we also
have $\pi'\in\lpw$.  The map $\pi\mapsto \pi'$ is a parity-reversing
involution (i.e., exactly one of $\pi$ and $\pi'$ is an even permutation)
on $\lpw$, and the proof follows. $\Box$

The above proof illustrates what will be our basic technique for
showing that a poset $P$ is sign-balanced, viz., giving a bijection
$\sigma:\lpw \rightarrow \lpw$ such that $\pi$ and $\sigma(\pi)$ have
opposite parity for all $\pi\in\lpw$. Equivalently, we are giving a
parity-reversing bijection $\vartheta:\ep\rightarrow\ep$.

In 1992 Ruskey \cite[{\S}5, item~6]{ruskey2} conjectured as to when
the product $\bm{m}\times \bm{n}$ of two chains of cardinalities $m$
and $n$ is sign-balanced, viz., $m,n>1$ and $m\equiv
n\,(\mathrm{mod}\,2)$. Ruskey proved this when $m$ and $n$ are both
even by giving a simple parity-reversing involution, which we
generalize in Proposition~\ref{prop:alch} and Corollary~\ref{cor:dom}.
Ruskey's conjecture for $m$ and $n$ odd was proved by D. White
\cite{white}, who also computed the ``imbalance'' between even and odd
linear extensions in the case when exactly one of $m$ and $n$ is even
(stated here as Theorem~\ref{thm:white}).
None of our theorems below apply to the case when $m$ and $n$ are both
odd. Ruskey \cite[{\S}5, item~5]{ruskey2} also asked what order ideals
$I$ (defined below) of $\bm{m}\times \bm{n}$ are sign-balanced. Such
order ideals correspond to integer partitions $\lambda$ and will be
denoted $P_\lambda$; the linear extensions of $P_\lambda$ are
equivalent to standard Young tableaux (SYT) of shape $\lambda$. White
\cite{white} also determined some additional $\lambda$ for which
$P_\lambda$ is sign-balanced, and our results below will give some
further examples.  In Sections~\ref{sec:maj} and \ref{sec:hl} we
consider some analogous questions for the parity of the major index of
a linear extension of a poset $P$.

Given $\pi=a_1a_2\cdots a_n\in\lpw$, let inv$(f)$ denote the number of
\emph{inversions} 
of $\pi$, i.e.,
  $$ \mathrm{inv}(\pi) = \#\{(i,j)\st i<j,\ a_i>a_j\}. $$
Let 
  \beq I_{P,\omega}(q) = \sum_{\pi\in\lpw}q^{\mathrm{inv}(f)},
   \label{eq:ipq} \eeq 
the generating function for linear extensions of $\ppw$ by number of
inversions. Since $f$
is an even linear extension if and only if inv$(f)$ is an even
integer, we see that $P$ is sign-balanced if and only if
$I_{P,\omega}(-1)=0$. In general $I_{P,\omega}(q)$ seems difficult to
understand, even when $P$ is known to be sign-balanced.

I am grateful to Marc van Leeuwen for his many helpful suggestions
regarding Section~\ref{sec:ptn}.

\section{Promotion and evacuation.}  
\begin{figure}
\centerline{\psfig{figure=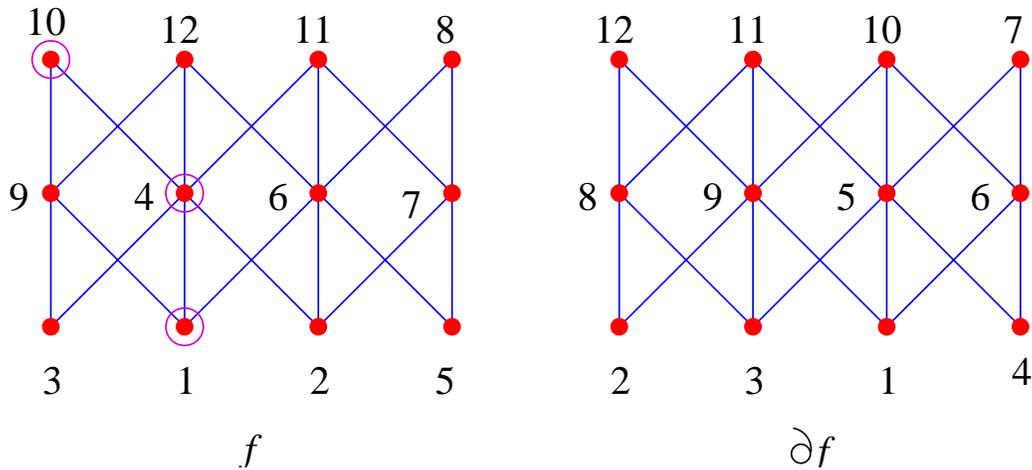}}
\caption{The promotion operator $\partial$}
\label{fig1}
\end{figure}

Promotion and evacuation are certain bijections on the set $\ep$ of
linear extensions of a finite poset $P$. They were originally defined
by M.-P.\ Sch\"utzenberger \cite{schut} and have subsequently arisen
is many different situations (e.g.,
\cite[{\S}5]{g-e}\cite[{\S}8]{haiman}\cite[{\S}4]{haiman2}\cite[{\S3}]{leeu}).
To be precise, the original definitions of promotion and evacuation
require an insignificant reindexing to become bijections. We will
incorporate this reindexing into our definition. Let $f:P\rightarrow
[n]$ be a linear extension of the poset $P$. Define a maximal chain
$u_0<u_1<\cdots<u_\ell$ of $P$, called the \emph{promotion chain} of $f$,
as follows. Let $u_0=f^{-1}(1)$. Once $u_i$ 
is defined let $u_{i+1}$ be that element $u$ covering $u_i$ (i.e.,
$u_i<u_{i+1}$ and no $s\in P$ satisfies $u_i<s<u_{i+1}$) for which
$f(u)$ is minimal. Continue until reaching a maximal element $u_\ell$ of
$P$. Now define the \emph{promotion} $g=\partial f$ of $f$ as follows.
If $t\neq u_i$ for any $i$, then set $g(t)=f(t)-1$. If $1\leq i\leq
k-1$, then set $g(u_i) = f(u_{i+1})-1$. Finally set
$g(u_\ell)=n$. Figure~\ref{fig1} gives an example, with the elements
in the promotion chain of $f$ circled. (The vertex labels in
Figure~\ref{fig1} are the values of a linear extension and are
unrelated to the (irrelevant) labeling $\omega$.) It is easy to see
that $\partial f\in\ep$ and that the map $\partial:\ep\rightarrow\ep$
is a bijection. 

\begin{lemma} \label{lemma1}
\emph{Let $P$ be an $n$-element poset. Then the promotion operator
$\partial: \ep\rightarrow \ep$ is parity-reversing if
and only if the length $\ell$ (or cardinality $\ell+1$) of every
maximal chain of $P$ satisfies $n\equiv \ell\modt$. Similarly,
$\partial$ is parity-preserving if and only if the length $\ell$ of
every maximal chain of $P$ satisfies $n\equiv \ell+1\modt$.}
\end{lemma}

\textbf{Proof.} Let $f\in\ep$, and let $u_0<u_1<\dots<u_\ell$ be the
promotion chain of $f$. Then $(\partial f)f^{-1}$ is a product of two
cycles, viz., 
  $$ (\partial f)f^{-1}= (n,n-1,\dots,1)(b_0,b_1,\dots,b_\ell), $$
where $b_i=f(u_i)$. This permutation is odd if and only if
$n\equiv\ell\,(\mathrm{mod}\,2)$, and the proof follows since every
maximal chain of $P$ is the promotion chain of some linear
extension. $\ \Box$   

\begin{corollary} \label{cor:sbmc}
\emph{Let $P$ be an $n$-element poset, and suppose that the length
$\ell$ of every maximal chain of $P$ satisfies $n\equiv
\ell\,(\mathrm{mod}\,2)$. Then $P$ is sign-balanced.}
\end{corollary}

\textbf{Proof.} By the previous lemma, $\partial$ is
parity-reversing. Since it is also a bijection,
$\ep$ must contain the same number of even linear extensions as odd
linear extensions. $\ \Box$

\begin{figure}
\centerline{\psfig{figure=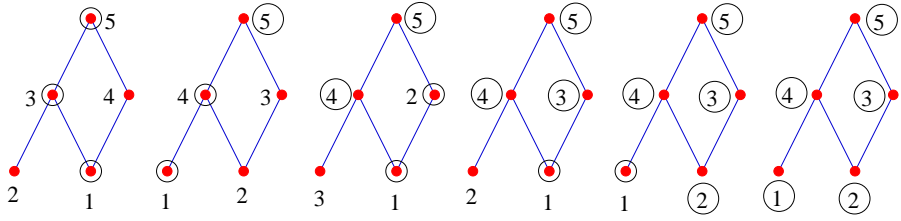}}
\caption{The evacuation operator evac.}
\label{fig2}
\end{figure}

We now consider a variant of promotion known as evacuation. For any
linear extension $g$ of an $m$-element poset $Q$, let
$u_0<u_1<\cdots<u_\ell$ be the promotion chain of $g$, so $\partial
g(u_\ell)=m$. Define $\rho_g(Q)=Q-\{u_\ell\}$. The restriction of
$\partial g$ to $\rho_g(Q)$, which we also denote by $\partial g$, is
a linear extension of $\rho_g(Q)$. Let 
  $$ \mu_{g,k}(Q)=\rho_{\partial^kg}\,\rho_{\partial^{k-1}g}\cdots
    \rho_{\partial g}\,\rho_g(Q). $$
Now let $\#P=n$ and define the \emph{evacuation} evac$(f)$
of $f$ to be the linear extension of $P$ whose value at the unique
element of $\mu_{g,k-1}(P)-\mu_{g,k}(P)$ is $n-k+1$, for $1\leq k\leq
n$. Figure~\ref{fig2} gives an example of evac$(f)$, where we circle
the values of evac$(f)$ as soon as they are determined. A remarkable
theorem of Sch\"utzenberger \cite{schut} asserts that evac is an
involution (and hence a bijection $\ep\rightarrow\ep$).

We say that the poset $P$ is \emph{consistent} if for all $t\in P$,
the lengths of all maximal chains of the principal order ideal
$\Lambda_t := \{s\in P\st s\leq t\}$ have the same parity. Let
$\nu(t)$ denote the length of the longest chain of $\Lambda_t$, and
set 
  $$ \Gamma(P) = \sum_{t\in P}\nu(t). $$ 
We also say that a permutation $\sigma$ of a finite set has
\emph{parity} $k\in\zz$ if either $\sigma$ and $k$ are both even or
$\sigma$ and $k$ are both odd. Equivalently, inv$(\sigma)\equiv
k\,(\mathrm{mod}\,2)$. 

\begin{proposition} \label{prop:conev}
\emph{Suppose that $P$ is consistent. Then} evac$:\ep \rightarrow\ep$
\emph{is parity-preserving if ${n\choose 2}-\Gamma(P)$ is even, and
parity-reversing if ${n\choose 2}-\Gamma(P)$ is odd.}
\end{proposition}

\textbf{Proof.} The evacuation of a linear extension $f$ of an
$n$-element poset $P$ consists of $n$ promotions
$\delta_1,\dots,\delta_n$, where $\delta_i$ is applied to a certain
subposet $P_{i-1}$ of $P$ with $n-i+1$ elements. Let $f_i$ be the linear
extension of $P$ whose restriction to $P_i$ agrees with $\delta_i
\delta_{i-1}\cdots \delta_1$, and whose value at the unique element
of $P_{j-1}-P_j$ for $j\leq i$ is $n-i+1$. Thus $f_0=f$ and $f_n=
\mathrm{evac}(f)$. (Figure~\ref{fig2} gives an example of the sequence
$f_0,\dots,f_5$.) Let $u_i$ be the end (top) of
the promotion chain for the promotion $\delta_i$. Thus
$\{u_1,u_2,\dots,u_n\}=P$.  Lemma~\ref{lemma1} shows that if $P$ is
consistent, then $f_if_{i-1}^{-1}$ has parity
$n-i+1-(\nu(u_i)+1)$. Hence the parity of evac$(f)f^{-1}$ is given by
  $$ \sum_{i=1}^n (n-i-\nu(u_i)) ={n\choose 2}-\sum_{t\in P}\nu(P)
   = {n\choose 2}-\Gamma(P), $$
from which the proof follows. $\ \Box$

\begin{corollary} \label{cor:cons}
\emph{Suppose that $P$ is consistent and ${n\choose 2}-\Gamma(P)$ is
  odd. Then $P$ is sign-balanced.} 
\end{corollary}

\textsc{Note.} In \cite[pp.\ 50--51]{rs:thesis}\cite[Cor.\ 
19.5]{rs:mem} it was shown using the theory of $P$-partitions that the
number $e(P)$ of linear extensions of $P$ is even if $P$ is graded of
rank $\ell$ (i.e., every maximal chain of $P$ has length $\ell$) and
$n-\ell$ is even, and it was stated that it would be interesting to
give a direct proof. Our Corollary~\ref{cor:sbmc} gives a direct proof
of a stronger result. Similarly in
\cite[Cor.~4.6]{rs:thesis}\cite[Cor.~19.6]{rs:mem} it was stated (in
dual form) that if for all $t\in P$ all maximal chains of $\Lambda_t$
have the same length, and if ${n\choose 2}-\Gamma(P)$ is odd, then
$e(P)$ is even. Corollary~\ref{cor:cons} gives a direct proof of a
stronger result. 

\section{Partitions.} \label{sec:ptn}
In this section we apply our previous results and obtain some new
results for certain posets corresponding to (integer) partitions.
We first review some notation and terminology concerning
partitions. Further details may be found in \cite[Ch.\ 7]{ec2}. Let
$\lambda= (\lambda_1,\lambda_2,\dots)$ be a partition of $n$, denoted
$\lambda\vdash n$ or $|\lambda|=n$. Thus
$\lambda_1\geq\lambda_2\geq\cdots\geq 0$ and 
$\sum \lambda_i =n$. We can identify $\lambda$ with its \emph{diagram}
$\{(i,j)\in \pp\times \pp\st 1\leq j\leq\lambda_i\}$. Let $\mu$ be
another partition such that $\mu\subseteq\lambda$, i.e., $\mu_i\leq
\lambda_i$ for all $i$. Define the \emph{skew partition} or \emph{skew
  diagram} $\lambda/\mu$ by
  $$ \lambda/\mu = \{ (i,j)\in \pp\times \pp\st \mu_i+1\leq j\leq
    \lambda_i\}. $$ 
Write $|\lm|=n$ to denote that $|\lambda|-|\mu|=n$, i.e., $n$ is the
number of squares in the shape $\lm$, drawn as a Young diagram
\cite[p.\ 29]{ec1}. We can regard $\lm$ as a subposet of
$\pp\times\pp$ (with the usual coordinatewise ordering). We write
$P_\lm$ for this poset. As a set it is the same as $\lm$, but the
notation $P_\lm$ emphasizes that we are considering it to be a poset.
In this section we will only be concerned with ``ordinary'' shapes
$\lambda$, but in Section~\ref{sec:maj} skew shapes $\lambda/\mu$ will
arise as a special case of Proposition~\ref{prop:slabps}. 

The posets $P_\lambda$ are consistent for any $\lambda$, so we can ask
for which $P_\lambda$ is evacuation parity-reversing, i.e., ${n\choose
2}-\Gamma(P_\lambda)$ is odd. To this end, the \emph{content}
$c(i,j)$ of the cell $(i,j)$ is defined by $c(i,j)=j-i$
\cite[p.~373]{ec2}. Also let ${\cal O}(\mu)$ denote the number of odd
parts of the partition $\mu$. An \emph{order ideal} of a poset $P$ is
a subset $K\subseteq P$ such that if $t\in K$ and $s<t$, then $s\in
K$. Similarly a \emph{dual order ideal} or \emph{filter} of $P$ is a
subset $F\subseteq P$ such that if $s\in F$ and $t>s$, then $t\in F$.
If we successively remove two-element chains from $P_\lambda$ which
are dual order ideals of the poset from which they are removed, then
eventually we reach a poset core$_2(P_\lambda)$, called the
\emph{2-core} of $P_\lambda$, that contains no dual order ideals which
are two-element chains.  The 2-core is \emph{unique}, i.e.,
independent of the order in which the dual order ideals are removed,
and is given by $P_{\delta_k}$ for some $k\geq 1$, where $\delta_k$
denotes the ``staircase shape'' $(k-1,k-2,\dots,1)$. For further
information see \cite[Exer.~7.59]{ec2}. 

\begin{proposition}
Let $\lambda\vdash n$. The following numbers all have the same
 parity. 
 \be  \item[(a)] $\Gamma(P_\lambda)$
  \item[(b)] $\sum_{t\in P_\lambda}c(t)$
  \item[(c)] $\frac 12({\cal O}(\lambda) - {\cal
    O}(\lambda'))$
  \item[(d)] $\frac 12(n-{k\choose 2})$, where
    ${k\choose 2}=\#\mathrm{core}_2(P_\lambda)$
  \ee 
Hence if $a_\lambda$ denotes any of the above four numbers, then
evacuation is partity-reversing on $P_\lambda$ if and only if
${n\choose 2}-a_\lambda$ is odd.
\end{proposition} 

\textbf{Proof.} It is easy to see that if $t\in P_\lambda$, then
$\nu(t) \equiv c(t)\,(\mathrm{mod}\,2)$. Hence (a) and (b) have the
same parity. It is well-known and easy to see \cite[Exam.\ 3, p.\
11]{macd} that 
   $$ \sum_{t\in P_\lambda}c(t) = \sum {\lambda_i\choose 2} -\sum 
    {\lambda'_i\choose 2}. $$
Since $\sum \lambda_i=\sum \lambda'_i$, we have
  $$ \sum_{t\in P_\lambda}c(t) = \frac 12\left( \sum \lambda_i^2
   -\sum \left(\lambda'_i\right)^2\right). $$
Since $a^2\equiv 0,1\,(\mathrm{mod}\,4)$ depending on whether $a$ is
even or odd, we see that (b) and (c) have the same parity. If we
remove from $P_\lambda$ a 2-element dual order ideal which is also a
chain, then we remove exactly one element 
with an odd content. A 2-core is self-conjugate and hence has an even
content sum. Hence the number of odd contents of $P_\lambda$ is equal
to the number of dominos that must be removed from $P_\lambda$ in
order to reach core$_2(P_\lambda)$. It follows that (b) and (c)
have the same parity, completing the proof. $\ \Box$  

It can be shown \cite{rs:amm} that if $t(n)$ denotes the number of
partitions $\lambda\vdash n$ for which $a_\lambda$ is even, then
$t(n)= \frac 12(p(n)+f(n))$, where $p(n)$ denotes the total number of
partitions of $n$ and
  $$ \sum_{n\geq 0}f(n)x^n = \prod_{i\geq 1}\frac{1+x^{2i-1}}
   {(1-x^{4i})(1+x^{4i-2})^2}. $$
Hence the number $g(n)$ of partitions $\lambda\vdash n$ for which evac
is parity-reversing on $P_\lambda$ is given by
  $$ g(n) = \left\{ \begin{array}{rl}
   \frac 12(p(n)+f(n)), & \mathrm{if}\ {n\choose 2}\ \mathrm{is\ odd}\\[.1in]
   \frac 12(p(n)-f(n)), & \mathrm{if}\ {n\choose 2}\ \mathrm{is\ even}
   \end{array} \right. $$

We conclude this section with some applications of the theory of
domino tableaux. A \emph{standard domino tableau} (SDT) of shape
$\lambda\vdash 2n$ is a sequence 
  $$ \emptyset=\lambda^0\subset \lambda^1 \subset \cdots\subset
    \lambda^n=\lambda $$
of partitions such that each skew shape $\lambda^i/\lambda^{i-1}$ is a
\emph{domino}, i.e., two squares with an edge in common. Each of these
dominos is either horizontal (two squares in the same row) or vertical
(two squares in the same column). 
Let $\dla$ denote the set of all SDT of shape $\lambda$.  
Given $D\in\dla$, define ev$(D)$ to be the number of vertical dominos
in even columns of $D$, where an \emph{even column} means the $2i$th
column for some $i\in\pp$. For the remainder of this section,
fix the labeling $\omega$ of $P_\lambda$ to be the usual ``reading
order,'' i.e., the first row of $\lambda$ is labelled
$1,2,\dots,\lambda_1$; the second row is labelled $\lambda_1+1,
\lambda_1+2,\dots, \lambda_1+\lambda_2$, etc. We write $I_\lambda(q)$
for $I_{P_\lambda,\omega}(q)$ and set $I_\lambda=I_\lambda(-1)$, the
\emph{imbalance} of the partition $\lambda$. It is shown in
\cite[Thm.\ 12]{white} (by analyzing the formula that results from
setting $q=-1$ in (\ref{eq:oiip})) that 
  $$ I_\lambda =  \sum_{D\in\dla} (-1)^{\mathrm{ev}(D)}. $$
\indent Let $\lambda\vdash n$. Lascoux, Leclerc and Thibon
\cite[(27)]{l-l-t} define a certain class of symmetric functions
$\tilde{G}^{(k)}_\lambda(x;q)$ (defined earlier by Carr\'e and Leclerc
\cite{c-l} for the special case $k=2$ and $\lambda=2\mu$). We will
only be concerned with the case $k=2$ and $q=-1$, for which we write
$G_\lambda = \tilde{G}^{(2)}_\lambda(x;-1)$. The symmetric function
$G_\lambda$ vanishes unless core$_2(\lambda)=\emptyset$, so we may
assume $n=2m$. If core$_2(\lambda)=\emptyset$, then $G_\lambda$ is
homogeneous of degree $m=n/2$. We will not define it here but only recall 
the properties relevant to us. The connection with the imbalance
$I_\lambda$ is provided by the formula (immediate from the definition
of $G_\lambda$ in \cite{l-l-t} together with \cite[Thm. 12]{white})
  \beq [x_1\cdots x_m]G_\lambda = (-1)^{r(\lambda)}I_\lambda,
    \label{eq:sfg} \eeq
where $[x_1\cdots x_m]F$ denotes the coefficient of $x_1\cdots x_m$ in
the symmetric function $F$, and $r(\lambda)$ is the maximum number of
vertical dominos that can appear in even columns of a domino tableau
of shape $\lambda$.  Also define $d(\lambda)$ to be the maximum number
of disjoint vertical dominos that can appear in the diagram of
$\lambda$, i.e., 
  $$ d(\lambda) = \sum_i\left\lfloor \frac
    12\lambda'_{2i}\right\rfloor. $$ 
Note that $d(\lambda)\geq r(\lambda)$, but equality need not hold in
general. For instance, $d(4,3,1)=1$, $r(4,3,1)=0$. However, we do have
$d(2\mu)= r(2\mu)$ for any partition $\mu$. Let us also note
that our $r(\lambda)$ is denoted $d(\lambda)$ in \cite{white} and is
defined only for $\lambda$ with an empty 2-core. 

\begin{theorem} \label{thm:kcor}
(a) \emph{We have}
  $$ \sum_{\mu\vdash m} I_{2\mu} = 1 $$
\emph{for all $m\geq 1$.}\\
\indent (b) \emph{Let $v(\lambda)$ denote the maximum number of
  disjoint vertical dominos that fit in the shape
  $\lambda$. Equivalently,} 
  $$ v(\lambda) = \sum_{i\geq 1} \left\lfloor \frac 12\lambda_i^\prime
     \right\rfloor. $$
\emph{Then} 
  $$ \sum_{\lambda\vdash 2m} (-1)^{v(\lambda)}I_\lambda^2 =0. $$ 
\end{theorem}

\textbf{Proof.} (a) Barbasch and Vogan \cite{b-v} and Garfinkle \cite{gar} 
define a bijection between elements $\pi$ of the hyperoctahedral group
$B_m$, regarded as signed permutations of $1,2,\dots,m$, and pairs
$(P,Q)$ of SDT of the same shape $\lambda\vdash 2m$.  (See
\cite[p.~25]{leeuwen} for further information.) A crucial property of
this bijection, stated implicitly without proof in \cite{k-l-l-t} and
proved by Shimozono and White \cite[Thm.\ 30]{s-w}, asserts that
  \beq \mathrm{tc}(\pi) = \frac 12(v(P)+v(Q)), \label{eq:srsk} \eeq
where tc$(\pi)$ denotes the number of minus signs in $\pi$ and $v(R)$
denotes the number of vertical dominos in the SDT $R$. 

Carr\'e and Leclerc \cite[Def.\ 9.1]{c-l} define a
symmetric function $H_\mu(x;q)$ which satisfies
$H_\mu(x,-1)=(-1)^{v(\mu)} G_{2\mu}$. In \cite[Thm.\ 1]{k-l-l-t} is
stated the identity 
  \beq \sum_\mu H_\mu(x;q)=\prod_i \frac{1}{1-x_i} \prod_{i<j}
    \frac{1}{1-x_ix_j}\prod_{i\geq j}\frac{1}{1-qx_ix_j}. 
    \label{eq:symqc} \eeq
The proof of (\ref{eq:symqc}) in \cite{k-l-l-t} is incomplete, since
it depends on a semistandard version of the $P=Q$ case of
(\ref{eq:srsk}) (easily deduced from (\ref{eq:srsk})), which had not
yet been proved. The proof of (\ref{eq:srsk}) in \cite{s-w} therefore
completes the proof of (\ref{eq:symqc}). A generalization of
(\ref{eq:symqc}) was later given by Lam \cite[Thm.~28]{lam}.

Setting $q=-1$ in (\ref{eq:symqc}) gives
  $$ \sum_{\mu} (-1)^{v(\mu)} G_{2\mu} = \prod_i
  \frac{1}{(1-x_i)(1+x_i^2)} \prod_{i<j} \frac{1}{1-x_i^2 x_j^2}. $$
Taking the coefficient of $x_1\cdots x_m$ on both sides and using
(\ref{eq:sfg}) together with $v(\mu)=d(2\mu)=r(2\mu)$ completes the
proof. 

(b) It is easy to see that for any SDT $D$ we have 
  $$ v(D) = v(\lambda)-2d(\lambda)+
      2\mathrm{ev}(D). $$
Thus by (\ref{eq:srsk}) we have
  \beas 0 & = & \sum_{\pi\in B_m}(-1)^{\mathrm{tc}(\pi)}\\
   & = & \sum_{P,Q}(-1)^{\frac 12(v(P)+v(Q))}\\ & = &
   \sum_{\lambda\vdash 2m}\left(\sum_{D\in\dla}(-1)^{\frac 12 v(D)}
   \right)^2\\ & = & \sum_{\lambda\vdash 2m} (-1)^{v(\lambda)}
   \left(\sum_{D\in\dla}(-1)^{\mathrm{ev}(D)}\right)^2\\ 
   & = & \sum_{\lambda\vdash 2m} (-1)^{v(\lambda)}I_\lambda^2.
   \ \ \Box \eeas
In the same spirit as Theorem~\ref{thm:kcor} we have the following
conjecture.

  \begin{conjecture} \label{conj:sytimb} \hspace{-.2in}\footnote{A
  combinatorial proof of (a) was found by Thomas Lam \cite{lam} after
  this paper was written. Later a combinatorial proof of both (a) and 
  (b) was given by Jonas Sj\"ostrand \cite{sjo}. Sj\"ostrand's main
  result \cite[Thm.~2.3]{sjo} leads to further identities, such as
  $\sum_{\mu\vdash n}q^{v(\mu)}I_{2\mu}=1$, thereby generalizing our
  Theorem~\ref{thm:kcor}(a).}
(a) For all $n\geq 0$ we have 
  \beq \sum_{\lambda\vdash n} q^{v(\lambda)}t^{d(\lambda)}
   x^{v(\lambda')}y^{d(\lambda')} I_\lambda 
    = (q+x)^{\lfloor n/2\rfloor}. \label{eq:opq} \eeq
(b) If $n\not\equiv 1\,(\mathrm{mod}\,4)$, then
  $$ \sum_{\lambda\vdash n}
   (-1)^{v(\lambda)}t^{d(\lambda)}I_\lambda^2=0. $$  
  \end{conjecture}

\begin{figure}
\centerline{\psfig{figure=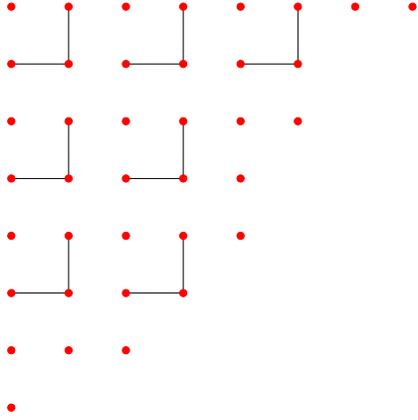}}
\caption{$d(86655431)=d(86655431')$}
\label{fig:horver}
\end{figure}

It is easy to see that $d(\lambda)=d(\lambda')$ for all
$\lambda$. (E.g., consider the horizontal and vertical line segments
in Figure~\ref{fig:horver}.) Hence the variable $y$ is superfluous in
equation (\ref{eq:opq}), but we have included it for the sake of
symmetry. In particular, if $F_n(q,t,x,y)$ denotes the left-hand side
of (\ref{eq:opq}) then 
  $$ F_n(q,0,x,y) = F_n(q,t,x,0) = F_n(q,0,x,0). $$
Note also that $d(\lambda)=0$ if and only $\lambda$ is a \emph{hook},
i.e., a partition of the form $(n-k,1^k)$. 

The case $t=0$ (or $y=0$, or $t=y=0$) of equation~(\ref{eq:opq})
follows from the following proposition, which in a sense
``explains'' where the right-hand side $(q+x)^{\lfloor n/2\rfloor}$
comes from.

\begin{proposition} \label{prop:hooksum}
\emph{For all $n\geq 0$ we have}
  \beq \sum_{\lambda=(n-k,1^k)} q^{v(\lambda)}
   x^{v(\lambda')} I_\lambda = (q+x)^{\lfloor n/2\rfloor}, 
   \label{eq:hsum} \eeq
\emph{where $\lambda$ ranges over all hooks $(n-k,1^k)$, $0\leq
  k\leq n-1$.}
\end{proposition}

\textbf{First proof.}  
Let $\lambda=(n-k,1^k)$. Let $\omega$ denote the
``reading order'' labeling of  $P_\lambda$ as above. The set ${\cal
L}_\pw$ consists of all permutations $1,a_2,\dots,a_m$, where
$a_2,\dots,a_m$ is a \emph{shuffle} of the permutations
$2,3,\dots,n-k$ and $n-k+1,n-k+2,\dots,n$. It follows e.g.\ from
\cite[Prop.~1.3.17]{ec1} that 
  $$ I_\lambda(q) = \qbc{n-1}{k}, $$
a $q$-binomial coefficient. 

Suppose first that $n=2m+1$. By \cite[Exer.~3.45(b)]{ec1},
  $$ \qbc{n-1}{k}_{q=-1}= \left\{ \begin{array}{rl} 
    \ds {m\choose j}, & k=2j\\[.2in] 0, & k=2j+1.
   \end{array} \right. $$
Note that if $\lambda=( n-2j,1^{2j})$, then
$v(\lambda)=j$ and $v(\lambda')=m-j$. Hence
  \beas \sum_{\lambda=( n-k,1^k)} q^{v(\lambda)}
    x^{v(\lambda')} I_\lambda & = & \sum_{j=0}^m q^j x^{m-j}{m\choose
      j}\\ & = & (q+x)^m, \eeas
as desired. The proof for $n$ even is similar and will be omitted. $\
\Box$ 

\textbf{Second proof.} Assume first that $n=2m$. We use an involution
argument analogous to the proof of Theorem~\ref{thm:p-r} or to
arguments in \cite[{\S}5]{white} and Section~\ref{sec:chains} of
this paper. Let $T$ be an SYT of shape $\lambda=
(n-k,1^k)$, which can be regarded as an element of ${\cal
L}_{P_\lambda,\omega}$. Let $i$ be the least positive integer (if it
exists) 
such that $2i-1$ and $2i$ appear in different rows and in different
columns of $T$. Let $T'$ denote the SYT obtained from $T$ by
transposing $2i-1$ and $2i$. Since multiplying by a transposition
changes the sign of a permutation, we have $(-1)^{\mathrm{inv}(T)} +
(-1)^{\mathrm{inv}(T')}=0$. The surviving SYT are obtained by first
placing $1,2$ in the same row or column, then $3,4$ in the same row or
column, etc. If $k=2j$ or $2j+1$, then the number of survivors is
easily seen to be ${m-1\choose j}$. Because the entries of $T$ come in
pairs $2i-1,2i$, the number of inversions of each surviving SYT is
even. Moreover, if $k=2j$ then $v(\lambda)=j$ and $v(\lambda')=m-j$,
while if $k=2j+1$ then $v(\lambda)=j+1$ and $v(\lambda')=m-1-j$. Hence
  \beas \sum_{\lambda=(n-k,1^k)} q^{v(\lambda)}
    x^{v(\lambda')}
    I_\lambda & = & \sum_{j=0}^{m-1}(q+x){m-1\choose j}q^j x^{m-1-j}\\
    & = & (q+x)^m, \eeas
as desired.

The proof is similar for $n=2m+1$. Let $i$ be the least positive
integer (if it exists) such that $2i$ and $2i+1$ (rather than $2i-1$
and $2i$) appear in different rows and in different columns of $T$.
There are now no survivors when $k=2j+1$ and ${m\choose j}$ survivors
when $k=2j$. Other details of the proof remain the same, so we get
  \beas \sum_{\lambda=(n-k,1^k)} q^{v(\lambda)}
    x^{v(\lambda')}
    I_\lambda & = & \sum_{j=0}^m {m-1\choose j}q^j x^{m-j}\\
    & = & (q+x)^m, \eeas 
completing the proof. $\ \Box$



There are some additional properties of the symmetric functions
$G_\lambda$ that yield information about $I_\lambda$. For instance,
there is a product formula in \cite[Thm.\ 2]{k-l-l-t} for $\sum_{\mu}
G_{2\mu\cup 2\mu}$, where $\mu$ ranges over all partitions and
  $$ 2\mu\cup 2\mu=(2\mu_1,2\mu_1,2\mu_2,2\mu_2, \dots), $$
which implies that $\sum_{\mu\vdash n}I_{2\mu\cup 2\mu}=0$. In
fact, in \cite[Cor.\ 9.2]{c-l} it is shown that $G_{2\mu\cup 2\mu}(x)=\pm 
s_\mu(x_1^2, x_2^2,\dots)$, from which it follows easily that in fact
$I_{2\mu\cup 2\mu}=0$. However, this result is just a special case of
Corollary~\ref{cor:sbmc} and of Proposition~\ref{prop:conev}, so we
obtain nothing new.

Also relevant to us is an expansion of $G_\lambda$ into Schur
functions due to Shimozono (see \cite[Thm.\ 18]{white}) for certain
shapes $\lambda$, namely, those whose 2-quotient (in the sense e.g.\ 
of \cite[Exam.~I.1.8]{macd}) is a pair of rectangles. This expansion
was used by White \cite[Cor.\ 20]{white} to evaluate $I_\lambda$ for
such shapes. White \cite[{\S}8]{white} also gives a combinatorial
proof, based on a sign-reversing involution, in the special case that
$\lambda$ itself is a rectangle. We simply state here White's result
for rectangles.

\begin{theorem} \label{thm:white}
\emph{Let $\lambda$ be an $m\times n$ rectangle. Then}
  $$ I_\lambda = \left\{ \begin{array}{rl} 1, & \mathrm{if}\ m=1\ 
    \mathrm{or}\ n=1\\
      0, & \mathrm{if}\ m\equiv n\,(\mathrm{mod}\,2)\ \mathrm{and}\ 
      m,n>1\\
    \pm g^\mu, & m\not\equiv n\,(\mathrm{mod}\,2), \end{array}
   \right. $$
\emph{where $g^\mu$ denotes the number of shifted standard tableaux
(as defined e.g.\ in \cite[Exam.~III.8.12]{macd}) of shape}
  $$ \mu = \left( \frac{m+n-1}{2}, \frac{m+n-3}{2}, \cdots,
      \frac{|n-m|+3}{2},\frac{|n-m|+1}{2}\right). $$
\emph{(An explicit ``hook length formula'' for any $g^\mu$ appears
  e.g.\ in the reference just cited.)}
\end{theorem}

It is natural to ask whether Theorem~\ref{thm:white} can be
generalized to other partitions $\lambda$. In this regard, A. Eremenko
and A. Gabrielov (private communication) have made a remarkable
conjecture. Namely, if we fix the number $\ell$ of parts and parity of
each part of $\lambda$, then there are integers $c_1,\dots,c_k$ and
integer vectors $\gamma_1,\dots, \gamma_k\in\zz^\ell$ such that
  $$ I_\lambda = \sum_{i=1}^k c_i g^{\frac 12(\lambda+\gamma_i)}. $$
One defect of this conjecture is that the expression for $I_\lambda$
is not unique. We can insure uniqueness, however, by the additional
condition that all the vectors $\gamma_i$ have coordinate sum 0 when
$|\lambda|$ is even and $-1$ when $|\lambda|$ is odd (where
$|\lambda|=\sum \lambda_i$). In this case, however, we need to define
properly $g^\mu$ when $\mu$ isn't a strictly decreasing sequence of
nonnegative integers. See the discussion preceding
Conjecture~\ref{conj:egsf}.  For instance, we have
  \beas I_{(2a,2b,2c)} & = & g^{(a,b,c)}-g^{(a+1,b,c-1)}\\
   I_{(2a+1,2b,2c)} & = & g^{(a,b,c)}+g^{(a+1,b-1,c)}\\
   I_{(2a,2b+1,c)} & = & 0\\
   I_{(2a,2b,2c+1)} & = & -g^{(a+1,b-1,c)}-g^{(a+1,b,c-1)}\\
   I_{(2a+1,2b+1,2c)} & = & g^{(a+1,b,c)}+g^{(a+1,b+1,c-1)}\\
   I_{(2a+1,2b,2c+1)} & = & 0\\
   I_{(2a,2b+1,2c+1)} & = & g^{(a+1,b,c)}+g^{(a,b+1,c)}\\
   I_{(2a+1,2b+1,2c+1)} & = & g^{(a,b+1,c)}+g^{(a+1,b+1,c-1)}\\
   I_{(2a,2b,2c,2d)} & = & g^{(a,b,c,d)}-g^{(a+1,b,c-1,d)}-
        g^{(a+1,b+1,c-1,d-1)}-2g^{(a+1,b,c,d-1)}. \eeas
It is easy to see that $I_{(2a,2b+1,c)}=I_{(2a+1,2b,2c+1)}=0$, viz.,
the 2-cores of the 
partitions $(2a,2b+1,c)$ and $(2a+1,b,2c+1)$ have more than one
square. More generally, we have verified by induction the formulas
for $I_\mu$ when $\ell(\mu)\leq 3$.  

We have found a (conjectured) symmetric function generalization of the
Eremenko-Gabrielov conjecture. If $f(x)$ is any symmetric function,
define
  $$ f(x/x) = f(p_{2i-1}\rightarrow 2p_{2i-1},\  p_{2i}\rightarrow 0).
  $$
In other words, write $f(x)$ as a polynomial in the power sums $p_j$
and substitute $2p_{2i-1}$ for $p_{2i-1}$ and 0 for $p_{2i}$. In
$\lambda$-ring notation, $f(x/x)=f(X-X)$. Let $Q_\mu$ denote
Schur's shifted $Q$-function \cite[{\S}3.8]{macd}. The $Q_\mu$'s form
a basis for the ring $\qq[p_1,p_3,p_5,\dots]$. Hence $f(x/x)$ can be
written uniquely as a linear combination of $Q_\mu$'s.

We mentioned above that the symmetric function $G_\lambda$ was
originally defined only when core$_2(\lambda)=\emptyset$. We can
extend the definition to any $\lambda$ as follows. The original
definition has the form 
  \beq G_\lambda(x) = \sum_D (-1)^{\mathrm{cospin}(D)}x^D, 
    \label{eq:glax} \eeq
summed over all semistandard domino tableaux of shape $\lambda$, where
cospin$(\lambda)$ is a certain integer and $x^D$ a certain monomial
depending on $\lambda$. If $\#\mathrm{core}_2(\lambda)=1$, then define
$G_\lambda$ exactly as in (\ref{eq:glax}), except that we sum over all
semistandard domino tableaux of the skew shape $\lambda/1$. If
$\#\mathrm{core}_2(\lambda)>1$, then define $G_\lambda=0$. (In certain
contexts it would be better to define $G_\lambda$ by (\ref{eq:glax}),
summed over all semistandard domino tableaux of the skew shape
$\lambda/\mathrm{core}_2(\lambda)$, but this is not suitable for our
purposes.) Equation~(\ref{eq:sfg}) then continues to hold for any
$\lambda\vdash n$, where $m=\lfloor n/2\rfloor$. 

We also need to define $G_\mu(x/x)$ properly when $\mu$ is not a
strictly decreasing sequence of positive integers. The following
definition seems to be correct, but perhaps some modification is
necessary. Let $\mu=(\mu_1,\dots,\mu_k)\in \zz^k$. Trailing 0's are
irrelvant and can be ignored, so we may assume $\mu_k>0$. If $\mu$ is
not a sequence of distinct nonnegative integers, then
$G_\mu(x/x)=0$. Otherwise $G_\mu(x/x) = \varepsilon_\mu
G_\lambda(x/x)$, where $\lambda$ is the decreasing rearrangement of
$\mu$ and $\varepsilon_\mu$ is the sign of the permutation that
converts $\mu$ to $\lambda$.

\begin{conjecture} \label{conj:egsf}
Fix the number $\ell$ of parts and parity of each part of the
partition $\lambda$. Then there are integers $c_1,\dots,c_k$ and
integer vectors $\gamma_1,\dots,\gamma_k\in\zz^\ell$ such that
  \beq (-1)^{r(\lambda)} G_\lambda(x/x) = \sum_{i=1}^k c_i Q_{\frac
         12(\lambda+\gamma_i)}(x). \label{eq:geconj} \eeq
\end{conjecture}

Let $\lambda\vdash 2n$ or $\lambda\vdash 2n+1$. Take the coefficient
of $x_1x_2\cdots x_n$ on both sides of (\ref{eq:geconj}). By
(\ref{eq:sfg}) the left-hand side becomes $2^nI_\lambda$. Moreover, if
$\mu\vdash m$ then the coefficient of $x_1\cdots x_m$ in $Q_\mu$ is
$2^m g^\mu$ \cite[(8.16)]{macd}. Hence Conjecture~\ref{conj:egsf}
specializes to the Eremenko-Gabrielov conjecture. At present we have
no conjecture for the values of the coefficients $c_i$. Here is a
short table (due to Eremenko and Gabrielov for $I_\lambda$; they have
extended this table to the case of four and five rows) of the
three-row case of Conjecture~\ref{conj:egsf}. For simplicity we write
$\pm$ for $(-1)^{r(\lambda)}$.
  \beas \pm G_{(2a,2b,2c)}(x/x) & = & Q_{(a,b,c)}(x) - 
     Q_{(a+1,b,c-1)}(x)\\ 
   \pm G_{(2a+1,2b,2c)}(x/x) & = & Q_{(a,b,c)}(x) +
     Q_{(a+1,b-1,c)}(x)\\
  \pm G_{(2a,2b+1,2c)}(x/x) & = & 0\\
  \pm G_{(2a,2b,2c+1)}(x/x) & = & -Q_{(a+1,b-1,c)}(x)
     -Q_{(a+1,b,c-1)}(x)\\
  \pm G_{(2a+1,2b+1,2c)}(x/x) & = & Q_{(a+1,b,c)}(x) +
      Q_{(a+1,b+1,c-1)}(x)\\
  \pm G_{(2a+1,2b,2c+1)}(x/x) & = & 0\\
  \pm G_{(2a,2b+1,2c+1)}(x/x) & = & Q_{(a+1,b,c)}(x)+
    Q_{(a,b+1,c)}(x)\\
  \pm G_{(2a+1,2b+1,2c+1)}(x/x) & = & Q_{(a,b+1,c)}(x)
    +Q_{(a+1,b+1,c-1)}(x). \eeas 
\indent We now discuss some general properties of the polynomial
$I_\lambda(q)$ and its value $I_\lambda(-1)$. Let $C(\lambda)$ denote
the set of \emph{corner squares} of $\lambda$, i.e., those squares of
the Young diagram of $\lambda$ whose removal still gives a Young
diagram. Equivalently, Pieri's formula \cite[Thm.~7.15.7]{ec2} implies
that
 \beq s_{\lambda/1}=\sum_{t\in C(\lambda)}
 s_{\lambda-t}. \label{eq:pieri} \eeq 
Let $f^\lambda$ denote the number of SYT of shape $\lambda$
\cite[Prop.~7.10.3]{ec2}, so
  \beq f^\lambda=\sum_{t\in C(\lambda)}f^{\lambda-t}. \label{eq:flrec}
  \eeq
Note that $I_\lambda(1)=f^\lambda$, so $I_\lambda(q)$ is a
(nonstandard) $q$-analogue of $f^\lambda$. The $q$-analogue of
equation (\ref{eq:flrec}) is the following result.

\begin{proposition} \label{prop:ilqrec}
\emph{We have}
 $$ I_\lambda(q) = \sum_{t\in C(\lambda)} q^{b_\lambda(t)}
    I_{\lambda-t}(q), $$
\emph{where $b_\lambda(t)$ denotes the number of squares
in the diagram of $\lambda$ in a lower row than that of $t$.}
\end{proposition}

\textbf{Proof.} We have by definition
  $$ I_\lambda(q) = \sum_T q^{\mathrm{inv}(\pi(T))}, $$
where $T$ ranges over all SYT of shape $\lambda$ and $\pi(T)$ is the
permutation obtained by reading the entries of $T$ in the usual
reading order, i.e., left-to-right and top-to-bottom when $T$ is
written in ``English notation'' as in
\cite{macd}\cite{ec1}\cite{ec2}. Suppose 
$\lambda\vdash n$. If $T$ is an SYT of shape $\lambda$, then the 
square $t$ occupied by $n$ is a corner square. The number of
inversions $(i,j)$ of $\pi(T)=a_1\cdots a_m$ such that $a_i=n$ is
equal to $b_\lambda(t)$, and the proof follows. $\ \Box$

Now let $D_1$ denote the linear operator on symmetric functions
defined by $D_1(s_\lambda)=s_{\lambda/1}$. We then have the commutation
relation \cite[Exercise~7.24(a)]{ec2}
  \beq D_1 s_1 - s_1 D_1 = I, \label{eq:dssd} \eeq
the identity operator. This leads to many enumerative consequences,
discussed in \cite{rs:dp}. There is an analogue of (\ref{eq:dssd})
related to $I_\lambda$, though we don't know of any
applications. Define a linear operator $D(q)$ on symmetric functions
by 
  $$ D(q)s_\lambda =\sum_{t\in C(\lambda)}
  q^{b_\lambda(t)}s_{\lambda-t}. $$
Let $U(q)$ denote the adjoint of $D(q)$ with respect to the basis
$\{s_\lambda\}$ of Schur functions, so
  $$ U(q)s_\mu=\sum_t q^{b_{\mu+t}(t)}s_{\mu+t}, $$
where $t$ ranges over all boxes that we can add to the diagram of
$\mu$ to get the diagram of a partition $\mu+t$ (for which necessarily
$t\in C(\mu+t)$). Note that $U(1)=s_1$ (i.e., multiplication by $s_1$)
and $D(1) = D_1$ as defined above. It
follows from Proposition~\ref{prop:ilqrec} that
  $$ U(q)^n\cdot 1 = \sum_{\lambda\vdash n} I_\lambda(q)s_\lambda, $$
where $U(q)^n\cdot 1$ denotes $U(q)^n$ acting on the symmetric
function $1=s_{\emptyset}$.  Write $U=U(-1)$ and $D=D(-1)$. Let $A$ be
the linear operator on symmetric functions given by $As_\lambda =
(2k(\lambda)+1)s_\lambda$, where $k(\lambda)=\#C(\lambda)$, the
number of corner boxes of $\lambda$.

\begin{proposition}
\emph{We have $DU+UD=A$.}
\end{proposition}

\textbf{Proof.} The proof is basically a brute force
computation. Write $\bar{\lambda}_i = \lambda_i+\lambda_{i+1}
+\cdots$. Suppose $\mu$ is obtained from $\lambda$ by adding a box
in row $r-1$ and deleting a box in row $s-1$, where $r<s$. Then the
coefficient of $s_\mu$ in $(D(q)U(q)+U(q)D(q))s_\lambda$ is given by
  $$ \langle s_\mu,(D(q)U(q)+U(q)D(q))s_\lambda\rangle =
   q^{\bar{\lambda}_r}q^{\bar{\lambda}_s}+
   q^{\bar{\lambda}_s}q^{\bar{\lambda}_r-1}, $$
which vanishes when $q=-1$. Similarly if $r>s$ we get
  $$ \langle s_\mu,(D(q)U(q)+U(q)D(q))s_\lambda\rangle =
   q^{\bar{\lambda}_s}q^{\bar{\lambda}_r+1}+
   q^{\bar{\lambda}_r}q^{\bar{\lambda}_s}, $$
which again vanishes when $q=-1$. On the other hand, if $\lambda=\mu$
we have
  \beas \langle s_\lambda,(D(q)U(q)+U(q)D(q))s_\lambda\rangle & = & 
     (c(\lambda)+1)q^{2\bar{\lambda}_r} +
     c(\lambda)q^{2\bar{\lambda}_r}\\ & = &
     (2c(\lambda)+1)q^{2\bar{\lambda}_r}. \eeas
When $q=-1$ the right-hand side become $2c(\lambda)+1$, completing the
proof. $\ \Box$
    
\section{Chains of order ideals.} \label{sec:chains}
Suppose that $P$ is an $n$-element poset, and let $\alpha=
(\alpha_1,\dots,\alpha_k)$ be a composition of $n$, i.e.,
$\alpha_i\in\pp =\{1,2,\dots\}$ and $\sum \alpha_i=n$. Define an
$\alpha$-\emph{chain} of order ideals of $P$ to be a chain 
  \beq \emptyset=K_0\subset K_1\subset\cdots\subset K_k=P 
    \label{eq:ac} \eeq
of order ideals satisfying $\#(K_i-K_{i-1})=\alpha_i$ for $1\leq i\leq
k$. The following result is quite simple but has a number of
consequences.

\begin{proposition} \label{prop:alch}
\emph{Let $P$ be an $n$-element poset and $\alpha$ a fixed composition
of $n$. Suppose that for every $\alpha$-chain (\ref{eq:ac}) of order
ideals of $P$, at least one subposet $K_i-K_{i-1}$ is
sign-balanced. Then $P$ is sign-balanced.}
\end{proposition}

\textbf{Proof.} Let ${\cal C}$ be the $\alpha$-chain (\ref{eq:ac}). 
We say that a linear extension $f$ is
${\cal C}$-\emph{compatible} if 
 $$ K_1=f^{-1}(\{ 1,\dots,\alpha_1\}),\ \ K_2-K_1=
  f^{-1}(\{\alpha_1+1,\dots,\alpha_1+\alpha_2\}), $$
etc. Let inv$({\cal C})$ be the minimum number
of inversions of a ${\cal C}$-compatible linear extension. Clearly
  $$ \sum_f q^{\mathrm{inv}(f)} = q^{\mathrm{inv}({\cal C})}
   \prod_{i=1}^k I_{K_i-K_{i-1}}(q), $$
where the sum is over all ${\cal C}$-compatible $f$. Since every
linear extension is compatible with a unique $\alpha$-chain, there
follows
  \beq I_{P,\omega}(q) = \sum_{{\cal C}} q^{\mathrm{inv}({\cal C})}
   \prod_{i=1}^k I_{K_i-K_{i-1}}(q), \label{eq:oiip} \eeq
where ${\cal C}$ ranges over all $\alpha$-chains of order ideals of
$P$. The proof now follows by setting $q=-1$. $\ \Box$

Define a finite poset $P$ with $2m$ elements to be \emph{tilable by
dominos} if there is a chain $\emptyset=K_0\subset K_1\subset
\cdots\subset K_m=P$ of order ideals such that each subposet
$K_i-K_{i-1}$ is a two-element chain. Similarly, if $\#P=2m+1$ and
$1\leq j\leq m+1$ then we say that $P$ is \emph{$j$-tilable by
dominos} if there is a chain $\emptyset=K_0\subset K_1\subset
\cdots\subset K_{m+1}=P$ of order ideals such that $\#(K_i-K_{i-1})=2$
if $1\leq i\leq m+1$ and $i\neq j$ (so $\#(K_j-K_{j-1})=1$). Note that
being tilable by dominos is stronger than the existence of a partition
of $P$ into cover relations (or two element saturated chains). For
instance, the poset $P$ with cover relations $a<c, b<c, a<d, b<d$ can
be partitioned into the two cover relations $a<c$ and $b<d$, but $P$
is not tilable by dominos. When $n=2m$, we define a \emph{$P$-domino
tableau} to be a chain $\emptyset=K_0\subset K_1\subset \cdots\subset
K_m=P$ of order ideals such that $K_i-K_{i-1}$ is a two-element chain
for $1\leq i\leq m$. Similarly, when $n=2m+1$, we define a
(standard) \emph{$P$-domino tableau} to be a chain
$\emptyset=K_0\subset K_1\subset \cdots\subset K_{m+1}=P$ of order
ideals such that $K_i-K_{i-1}$ is a two-element chain for $1\leq i\leq
m$ (so that $K_{m+1}-K_m$ consists of a single point). Thus for
$\lambda\vdash 2n$, a $P_\lambda$-domino tableau coincides with our
earlier definition of an SDT of shape $\lambda$. 

\begin{corollary} \label{cor:dom}
{Let $\#P=2m$, and assume that $P$ is not tilable by dominos. Then $P$
  is sign-balanced. Similarly if $\#P=2m+1\geq 3$ and  $P$ is not
  $j$-tilable by dominos for some $j$, then $P$ is sign-balanced.}
\end{corollary}

\textbf{Proof.} Let $\alpha=(2,2,\dots,2)$ ($m$ 2's). If $\#P=2m$ and
$P$ is not tilable by dominos, then for any $\alpha$-chain (\ref{eq:ac})
there is an $i$ for which $K_i-K_{i-1}$ consists of two disjoint
points. Since a poset consisting of two disjoint points is
sign-balanced, it follows from Proposition~\ref{prop:alch} that $P$ is
sign-balanced. The argument is similar for $\#P=2m+1$. $\ \Box$

Corollary~\ref{cor:dom} was proved in a special case (the product of
two chains with an even number of elements, with the $\hat{0}$ and
$\hat{1}$ removed), using essentially the same proof as we have given,
by Ruskey \cite[{\S}5, item~6]{ruskey2}.

Corollary~\ref{cor:dom} is particularly useful for the posets
$P_\lambda$. From this corollary and the definition of
core$_2(\lambda)$ we conclude the following.

\begin{corollary}
\emph{If} core$_2(P_\lambda)$ \emph{consists of more than one element,
  then $P_\lambda$ is sign-balanced.}
\end{corollary}

It follows from \cite[Exer.~7.59(e)]{ec2} that if $f(n)$ denotes
the number of partitions $\lambda\vdash n$
such that $\#$core$_2(\lambda)\leq 1$, then
  $$ \sum_{n\geq 0} f(n)x^n = \frac{1+x} {\prod_{i\geq
      1}(1-x^{2i})^2}. $$
Standard partition asymptotics (e.g., \cite[Thm.\ 6.2]{andrews}) shows 
that 
  $$ f(n) \sim \frac{C}{n^{5/4}}\exp\left( \pi\sqrt{2n/3}\right) $$
for some $C>0$.
Since the total number $p(n)$ of partitions of $n$ satisfies 
  $$ p(n) \sim \frac{C'}{n}\exp\left( \pi\sqrt{2n/3}\right), $$
it follows that $\lim_{n\geq 0} f(n)/p(n)=0$. Hence as $n\rightarrow
\infty$, $P_\lambda$ is sign-balanced for almost all 
$\lambda\vdash n$. 

\section{Maj-balanced posets.} \label{sec:maj}
If $\pi=a_1 a_2\cdots a_m$ is a permutation of $[n]$, then the
\emph{descent set} $D(\pi)$ of $\pi$ is defined as
  $$ D(\pi) = \{ i\st a_i>a_{i+1}\}. $$
An element of $D(\pi)$ is called a \emph{descent} of $\pi$, and
\emph{major index} maj$(\pi)$ is defined as     
  $$ \maj(\pi) = \sum_{i\in D(\pi)} i. $$ 
The major index has many properties analogous to the number of
inversions, e.g., a classic theorem of MacMahon states that inv and
maj are equidistributed on the symmetric group $\sn$
\cite{foata}\cite{f-s}. Thus 
it is natural to try to find ``maj analogues'' of the results of the
preceding sections. In general, the major index of a linear extension
of a poset can be more tractable or less tractable than the number of
inversions. Thus, for example, in Theorem~\ref{thm:majdom} we are able
to completely characterize naturally labelled maj-balanced posets. An
analogous result for sign-balanced partitions seems very difficult. On
the other hand, since multiplying a permutation by a fixed permutation
has no definite effect on the parity of the major index, many of the
results for sign-balanced posets are false (Theorem~\ref{thm:p-r},
Lemma~\ref{lemma1}, Proposition~\ref{prop:conev}).

Let $f$ be a linear extension of the labelled poset $\ppw$, and
let $\pi=\pi(f)$ be the associated permutation of $[n]$. In analogy to
our definition of inv$(f)$, define maj$(f)=\mathrm{maj}(\pi)$ and 
  $$ W_{P,\omega}(q) = \sum_{f\in\ep}q^{\mathrm{maj}(f)} =
     \sum_{\pi\in\lpw}q^{\maj(\pi)}. $$
We say that $\ppw$ is \emph{maj-balanced} if $W_{P,\omega}(-1)=0$,
i.e., if the number of linear extensions of $P$ with even major index
equals the number with odd major index. Unlike the situation for
sign-balanced posets, the property of being maj-balanced can depend on
the labeling $\omega$. Thus an interesting special case is that of
\emph{natural labelings}, for which $\omega(s)<\omega(t)$ whenever
$s<t$ in $P$. We write $W_P(q)$ for $W_{P,\omega}(q)$ when $\omega$ is
natural. It is a basic consequence of the theory of $P$-partitions
\cite[Thm.~4.5.8]{ec1} that $W_P(q)$ does not depend on the choice of natural
labeling of $P$. 

\begin{figure}
\centerline{\psfig{figure=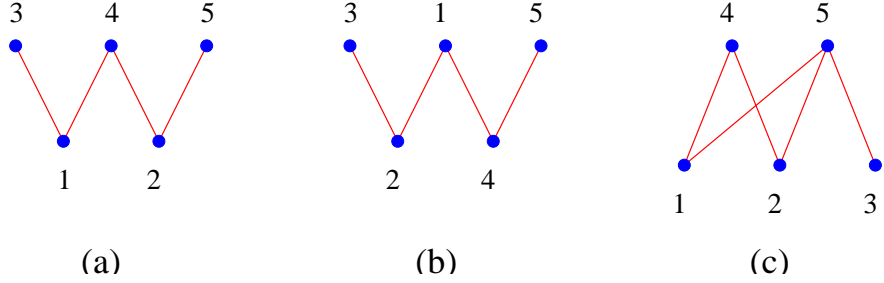}}
\caption{Some counterexamples}
\label{fig:cx}
\end{figure}

Figures~\ref{fig:cx}(a) and (b) show two different
labelings of a poset $P$. The first labeling (which is natural) is not
maj-balanced, while the second one is. Moreover, the dual poset $P^*$
to the poset $P$ in Figure~\ref{fig:cx}(b), whether naturally labelled or
labelled the same as $P$, is maj-balanced. Contrast that with the
trivial fact that the dual of a sign-balanced poset is sign-balanced.
As a further example of the contrast between sign and maj-balanced
posets, Figure~\ref{fig:cx}(c) shows a naturally labelled maj-balanced
poset $Q$. However, if we adjoin an element $\hat{0}$ below every
element of $Q$ and label it 0 (thus keeping the labeling natural) then
we get a poset which is no longer maj-balanced. On the other hand, it
is clear that such an operation has no effect on whether a poset is
sign-balanced. (In fact, it leaves $I_{Q,\omega}(q)$ unchanged.)

Corollary~\ref{cor:dom} carries over to the major index in the
following way. 

\begin{theorem} \label{thm:majdom}
(a) \emph{Let $P$ be naturally labelled. Then $W_P(-1)$ is equal to the
  number of $P$-domino tableaux. In particular, $P$ is maj-balanced if
  and only if there does not exist a $P$-domino tableau.}

(b) \emph{A labelled poset $(\pw)$ is maj-balanced if there does not
  exist a $P$-domino tableau.} 
\end{theorem}

\textbf{Proof.} 
(a) Let $\pi=a_1 \cdots a_m\in\lpw$. Let $i$ be the least number (if
it exists) for which $\pi'=a_1 \cdots a_{2i} a_{2i+2} a_{2i+1}
a_{2i+3}\cdots a_m\in \lpw$. Note that $(\pi')^\prime=\pi$. Now
exactly one of $\pi$ and $\pi'$ has a descent at $2i+1$.  The only
other differences in the descent sets of $\pi$ and $\pi'$ occur
(possibly) for the even numbers $2i$ and $2i+2$. Hence
$(-1)^{\mathrm{maj}(\pi)}+(-1)^{\mathrm{maj}(\pi')} =0$. The surviving
permutations $\sigma=b_1\cdots b_m$ in $\lpw$ are exactly those for
which $\emptyset\subset \{b_1,b_2\} \subset
\{b_1,\dots,b_4\}\subset\cdots$ is a $P$-domino tableau with
$\omega^{-1}(b_{2i-1})<\omega^{-1}(b_{2i})$ in $P$. (If $n$ is even,
then the $P$-domino tableau ends as $\{b_1,\dots,b_{n-2}\}\subset P$,
while if $n$ is odd it ends as $\{b_1,\dots,b_{n-1}\}\subset P$.)
Since $\omega$ is natural we have $b_{2i-1}<b_{2i}$ for all $i$, so
maj$(\sigma)$ is even. Hence $W_P(-1)$ is equal exactly to the number
of $P$-domino tableaux.
 
(b) Regardless of the labeling $\omega$, if there does not exist a
$P$-domino tableau then there will be no survivors in the argument of
(a), so $W_P(-1)=0$. $\ \Box$

\begin{figure}
\centerline{\psfig{figure=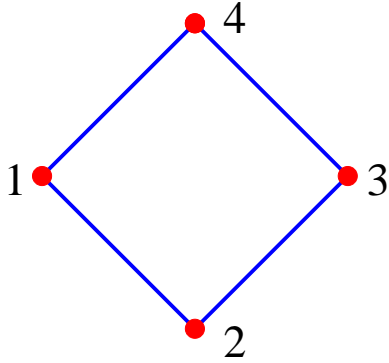}}
\caption{A maj-balanced labelled poset tilable by dominos}
\label{fig:notdt}
\end{figure}

The converse to Theorem~\ref{thm:majdom}(b) is false. The labelled
poset $(P,\omega)$ of Figure~\ref{fig:notdt} is tilable by dominos and
is maj-balanced. 

Given an $n$-element poset $P$ with dual $P^*$, set
$\Delta(P)=\Gamma(P^*)$. In
\cite[Thm~4.4]{rs:thesis}\cite[Prop.~18.4]{rs:mem}\cite[Thm.~4.5.2]{ec1}
it is shown that the following two conditions are equivalent:
  \be \item[(i)] For all $t\in P$, all maximal chains of the principal
  dual order ideal $V_t=\{s\in P\st s\geq t\}$ have the same length.
 \item[(ii)] $q^{{n\choose 2}-\Delta(P)}W_P(1/q)=W_P(q)$.
 \ee
It follows by setting $q=-1$ that if (i) holds and ${n\choose
2}-\Delta(P)$ is odd, then $P$ is
maj-balanced. Corollary~\ref{cor:cons} suggests in fact the following
stronger result. 

\begin{corollary} \label{cor:dcmb}
\emph{Suppose that $P$ is naturally labelled and dual consistent
(i.e., $P^*$ is consistent). If ${n\choose 2}-\Delta(P)$ is odd, then
$P$ is maj-balanced.}
\end{corollary}

\textbf{Proof.} By Theorem~\ref{thm:majdom} we need to show that there
does not exist a $P$-domino tableau. Given $t\in P$, let $\delta(t)$
denote the length of the longest chain of $V_t$, so $\Delta(P)=
\sum_{t\in P}\delta(t)$. First suppose that $n=2m$, and assume to the
contrary that $\emptyset=I_0\subset I_1\subset\cdots\subset I_m=P$ is
a $P$-domino tableau. If $s,t\in I_i-I_{i-1}$ then by dual consistency
$\delta(s)+\delta(t)\equiv 1\,(\mathrm{mod}\,2)$. Hence
$\Delta(P)\equiv m\,(\mathrm{mod}\,2)$, so
  $$ {n\choose 2}-\Delta(P) \equiv m(2m-1)-m\equiv
    0\,(\mathrm{mod}\,2), $$
a contradiction. 

Similarly if $n=2m+1$, then the existence of a $P$-domino tableau
implies $\Delta(P)\equiv m\,(\mathrm{mod}\,2)$, so
  $$ {n\choose 2}-\Delta(P) \equiv m(2m+1)-m\equiv
    0\,(\mathrm{mod}\,2), $$
again a contradiction. $\ \Box$

\begin{figure}
\centerline{\psfig{figure=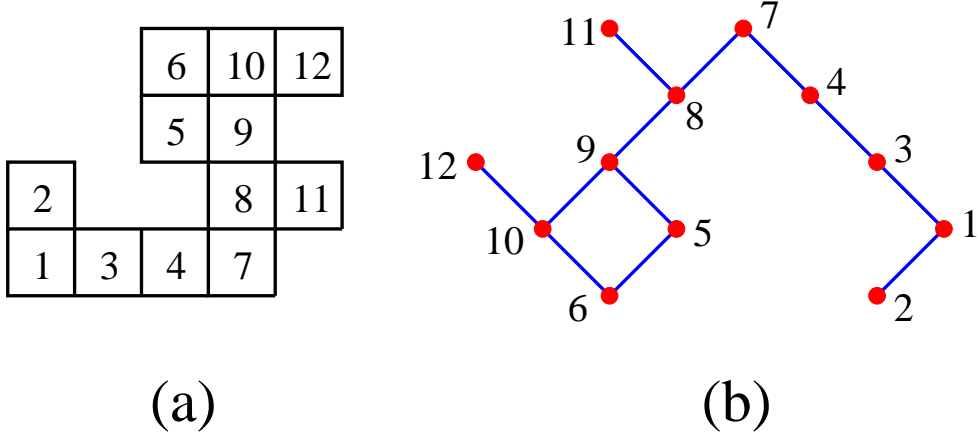}}
\caption{A set ${\cal S}$ of squares and the Schur labelled poset
  $P_\cs$} 
\label{fig:pposet}
\end{figure}

Now let ${\cal S}$ be a finite subset of solid unit squares with
integer vertices in $\rr\times\rr$ such that the set $|\cs|
=\bigcup_{S\in \cs}$ is simply-connected. For $S,T\in\cs$, define $S<T$
if the center vertices $(s_1,s_2)$ of $S$ and $(t_1,t_2)$ of $T$
satisfy either (a) $t_1=s_1$ and $t_2=s_2+1$ or (b) $t_1=s_1+1$ and
$t_2=s_2$. Regard $\cs$ as a poset, denoted $P_\cs$, under the
transitive (and reflexive) closure of the relation
$<$. Figure~\ref{fig:pposet} gives an example, where (a) shows $\cs$
as a set of squares and (b) as a poset. Note that the posets $P_\lm$
are a special case.

A \emph{Schur labelling} $\omega$ of $P_\cs$ is a labeling that
increases along rows and decreases along columns, as illustrated in
Figure~\ref{fig:pposet}. For the special case $P_\lm$, Schur labelings
play an important role in the expansion of skew Schur functions
$s_\lm$ in terms of quasisymmetric functions \cite[pp.~360--361]{ec2}.
Suppose that $\#P_\cs$ is even and that $P_\cs$ is tilable by dominos.
Then $\cs$ itself is tilable by dominos in the usual sense. It is
known (implicit, for instance, in \cite{thurston}, and more explicit
in \cite{chaboud}) that any two domino   
tilings of $\cs$ can be obtained from each other by ``$2\times 2$
flips,'' i.e., replacing two horizontal dominos in a $2\times 2$
square by two vertical dominos or \emph{vice versa}. It follows that
if $D$ is a domino tiling of $\cs$ with $v(D)$ vertical dominos, then
$(-1)^{v(D)}$ depends only on $\cs$. Set sgn$(\cs)=(-1)^{v(D)}$ for
any domino tiling of $\cs$.

\begin{proposition} \label{prop:slabps}
\emph{Let $\cs$ be as above, and let $\omega$ be a Schur labeling of
  $P_\cs$, where $\#P_\cs$ is even, say $\#P_\cs=n$. Then}
  sgn$(\cs)W_{P_\cs}(-1)$ 
  \emph{is the number of $P_\cs$-domino tableaux.} 
\end{proposition}

\textbf{Proof.} The proof parallels that of
Theorem~\ref{thm:majdom}. Define the involution $\pi\mapsto \pi'$ as
in the proof of Theorem~\ref{thm:majdom}. Each survivor $\sigma=b_1\cdots
b_m$ corresponds to a $P_\cs$-domino tableau $D$. We have
$b_{2i-1}>b_{2i}$ if and only if the domino labelled with $b_{2i-1}$
and $b_{2i}$ is vertical. As noted above,
$(-1)^{v(D)}=\mathrm{sgn}({\cal S})$, independent of $D$. Hence
$(-1)^{\maj(\sigma)}=\mathrm{sgn}(\sigma)$, and the proof follows as
in Theorem~\ref{thm:majdom}(a). $\ \Box$ 

A result analogous to Proposition~\ref{prop:slabps} holds for
$\#P_\cs$ odd (with essentially the same proof) provided $P_\cs$ has
a $\hat{0}$ or $\hat{1}$. The special case $P_\lm$ of
Proposition~\ref{prop:slabps} (and its analogue for $\#P_\cs$ odd) can
also be proved using the theory of symmetric functions, notably,
\cite[Prop.~7.19.11]{ec2} and the Murnaghan-Nakayama rule
(\cite[Cor.~7.17.5]{ec2}).

\section{Hook lengths.} \label{sec:hl}
In this section we briefly discuss a class of posets $P$ for which
$W_P(q)$, and sometimes even $I_{P,\omega}(q)$, can be explicitly
computed. For this class of posets we get a simple criterion for being
maj balanced and, if applicable, sign balanced.

Following \cite[p.\ 84]{rs:mem}, an $n$-element poset $P$ is called a
\emph{hook length poset} if there exist positive integers $h_1,\dots,
h_n$, the \emph{hook lengths} of $P$, such that 
  \beq W_P(q) = \frac{[n]!}{(1-q^{h_1})
    \cdots(1-q^{h_n})}, \label{eq:wpqhl} \eeq 
where $[n]!=(1-q)(1-q^2)\cdots(1-q^n)$.  It is easy to see that if $P$
is a hook length poset, then the multiset of hook lengths is unique.
In general, if $P$ is an ``interesting'' hook length poset, then each
element of $P$ should have a hook length associated to it in a
``natural'' combinatorial way.

\textsc{Note.} We could just as easily have extended our definition to
\emph{labelled} posets $(P,\omega)$, where now 
  $$ W_\pw(q) = \frac{q^c\,[n]!}{(1-q^{h_1})
      \cdots (1-q^{h_n})} $$
for some $c\in\nn$. However, little is known about the labelled
situation except when we can reduce it to the case of natural
labelings by subtracting certain constants from the values of
$\sigma$. 

The following result is an immediate consequence of equation
(\ref{eq:wpqhl}).

\begin{proposition} \label{prop:hlsb}
\emph{Suppose that $P$ is a hook length poset with hook lengths
  $h_1,\dots, h_n$. Then $P$ is maj-balanced if and only if the number
  of even hook lengths is less than $\lfloor n/2\rfloor$. If $P$ isn't
  maj-balanced, then the maj imbalance is given by}
  $$ W_P(-1) = \frac{\lfloor n/2\rfloor!}{\prod_{h_i\
      \mathrm{even}}(h_i/2)}. $$
\end{proposition}

It is natural to ask at this point what are the known hook length
posets. The strongest work in this area is due to Proctor
\cite{proc}\cite{proc2}. We won't state his remarkable results here,
but let us note that his \emph{$d$-complete} posets encompass all
known ``interesting'' examples of hook length posets. These include
forests (i.e., posets for which every element is covered by at most
one element) and the duals $P^\ast_\lambda$ of the posets $P_\lambda$
of Section~\ref{sec:ptn}.

Bj\"orner and Wachs \cite[Thm.\ 1.1]{b-w} settle the question of what
naturally labelled posets $(P,\omega)$ satisfy
 \beq I_{P,\omega}(q)=W_{P,\omega}(q). \label{eq:bw} \eeq
Namely, $P$ is a forest and $\omega$ is a postorder labeling.
Hence for postorder labelled forests, Proposition~\ref{prop:hlsb}
holds also for $I_{P,\omega}(-1)$. Bj\"orner and Wachs also obtain less
definitive results for arbitrary labelings, whose relevance to sign
and maj imbalance we omit.


\end{document}